\documentclass[11pt]{llncs}
\usepackage{epsfig}
\usepackage{graphicx}
\usepackage{color}
\usepackage{amsmath,amssymb}    
\usepackage{latexsym}
\usepackage{hhline}
\usepackage{algorithmic}
\usepackage{algorithm}
\spnewtheorem{fact}{Fact}{\bfseries}{\itshape}
  \baselineskip 8 mm

\newcommand{\ds}{\displaystyle}  
\newcommand{\comp}[1]{\overline{#1}}  

\newcommand{\cub}{\mbox{\textnormal{cub}}}
\newcommand{\boxi}{\mbox{\textnormal{box}}}
\newcommand{\ceil}[1]{\left\lceil #1 \right\rceil}

\newcommand{\gset}{\mathcal{X}}
\newcommand{\thresh}{\mathcal{T}}
\newcommand{\rlzr}{\mathcal{L}}
\newcommand{\ind}{\mathcal{I}}
\newcommand{\clq}{\mathcal{C}}

\begin{document}

\title{The Hardness of Approximating the Boxicity, Cubicity and Threshold
Dimension of a Graph}
\author{Abhijin Adiga\inst{1}, Diptendu Bhowmick\inst{1}, L. Sunil
Chandran\inst{1}}
\institute{Department of Computer Science and Automation, Indian Institute
of Science, Bangalore--560012, India. \\emails: abhijin@csa.iisc.ernet.in, diptendubhowmick@gmail.com, sunil@csa.iisc.ernet.in}
\date{}
\maketitle
{\renewcommand{\thefootnote}{}
\footnote{This work was supported by DST grant SR/S3/EECE/62/2006
and Infosys Fellowship.}}

\begin{abstract}
A $k$-dimensional box is the Cartesian product $R_1 \times R_2 \times
\cdots \times R_k$ where each $R_i$ is a closed interval on the
real line. The {\it boxicity} of a graph $G$, denoted as $\boxi(G)$,
is the minimum integer $k$ such that $G$ can be represented as the
intersection graph of a collection of $k$-dimensional boxes. A unit
cube in $k$-dimensional space or a $k$-cube is defined as the Cartesian
product $R_1 \times R_2 \times \cdots \times R_k$ where each $R_i$ is a
closed interval on the real line of the form $[a_i,a_i + 1]$. The {\it
cubicity} of $G$, denoted as $\cub(G)$, is the minimum integer $k$ such
that $G$ can be represented as the intersection graph of a collection
of $k$-cubes. The {\it threshold dimension} of a graph $G(V,E)$ is the
smallest integer $k$ such that $E$ can be covered by $k$ threshold
spanning subgraphs of $G$. In this paper we will show that there exists
no polynomial-time algorithm to approximate the threshold dimension of
a graph on $n$ vertices with a factor of $O(n^{0.5-\epsilon})$ for any $
\epsilon >0$, unless $NP=ZPP$. From this result we will show that there
exists no polynomial-time algorithm to approximate the boxicity and the
cubicity of a graph on $n$ vertices with factor $O(n^{0.5-\epsilon})$
for any $ \epsilon >0$, unless $NP=ZPP$. In fact all these hardness
results hold even for a highly structured class of graphs namely the
split graphs. We will also show that it is NP-complete to determine if
a given split graph has boxicity at most 3.
\end{abstract}

\section{Introduction}
Let $G(V,E)$ be a simple undirected finite graph with vertex set $V$ and
edge set $E$. A $d$-dimensional box is a Cartesian product $R_1\times
R_2\times\cdots\times R_d$ where each $R_i$ (for $1\leq i \leq d$) is
a closed interval of the form $[a_i,b_i]$ on the real line. A $k$-box
representation of $G$ is a mapping of the vertices of $G$
to $k$-boxes such that two vertices in $G$ are adjacent if and only
if their corresponding $k$-boxes have a non-empty intersection. The
\emph{boxicity} of a graph denoted $\boxi(G)$, is the minimum integer
 $k$ such that $G$ can be represented as the intersection graph of
$k$-dimensional boxes. A $d$-dimensional cube is a Cartesian product
$R_1\times R_2\times\cdots\times R_d$ where each $R_i$ (for $1\leq i \leq
d$) is a closed interval of the form $[a_i,a_i+1]$ on the real line. A
$k$-cube representation of a graph $G$ is a mapping of the vertices of
$G$ to $k$-cubes such that two vertices in $G$ are adjacent if and only
if their corresponding $k$-cubes have a non-empty intersection. The
\emph{cubicity} of $G$ is the minimum integer $k$ such that $G$ has a
$k$-cube representation.

The concept of boxicity was introduced by Roberts
 \cite{recentProgressesInCombRoberts}. Cozzens
\cite{phdThesisCozzens} showed that computing the boxicity of
a graph is NP-hard. This was later strengthened by Yannakakis
\cite{complexityPartialOrderDimnYannakakis} and finally by
Kratochv\'{\i}l \cite{specialPlanarSatisfiabilityProbNPKratochvil} who
showed that determining whether boxicity of a graph is at most two is
NP-complete. In \cite{complexityPartialOrderDimnYannakakis} Yannakakis
showed that it is NP-complete to determine whether the cubicity of a
given graph is at most 3. 

\subsection{Interval Graphs}
A graph $G$ is an \emph{interval graph} if and only if $G$ has an interval
representation: i.e. each vertex of $G$ can be associated with an interval
on the real line such that two intervals intersect if and only if the
corresponding vertices are adjacent. An interval graph $G$ is said to be a
\emph{unit interval graph} if and only if there is some interval representation of $G$ in which all the intervals are of the same length. Clearly, graphs with boxicity at most 1 are precisely the \emph{interval graphs} and the graphs with cubicity at most 1 are precisely the \emph{unit interval graphs}.

\subsection{Split Graphs}
A graph $G(V,E)$ is a split graph if its vertex set can be partitioned
into a clique and an independent set. We will denote the clique by $\clq(G)$
and independent set by $\ind(G)$. Note that this partition need not be
unique. But whenever we refer to $\clq(G)$, the set $V\setminus\clq(G)$
is an independent set and is denoted by $\ind(G)$.
Split graphs were first studied by F\"{o}ldes and Hammer in
\cite{splitGraphFoldes,aggregationInequalityChvatal}, and independently
introduced by Tyshkevich and Chernyak \cite{canonicalPartitionTyshkevich}.
For other characterizations and properties of split graphs one can refer
to Golumbic \cite{algGraphTheoryPerfectGraphsGolumbic}. 
\begin{fact}\label{fact:splitComp}
Complement of a split graph is a split graph.
\end{fact}

\begin{definition}
A split interval graph is a graph which is both a split graph and an
interval graph.
\end{definition}

\subsection{Threshold graphs and the Threshold Dimension Problem} 
A graph is a threshold graph if there is a
real number $S$ and a weight function $w: V \longrightarrow \mathbb{R}$
such that for any two vertices $u,v \in V(G)$, $(u,v)$ is an edge if
and only if $w(u) + w(v) \geq S$. We will use the following
property frequently in later sections.
\begin{fact} \label{fact:tgDefn}
A graph $G(V,E)$ is a threshold graph if and only if it is a split graph
and for every pair of vertices $u,v\in\ind(G)$, either $N(u)\subseteq N(v)$
or $N(v)\subseteq N(u)$. Equivalently, a threshold graph can be defined as
a split graph without an induced $P_4$ $($i.e. a path on 4 vertices$)$.
\end{fact}
Note that threshold graphs are interval graphs.
\begin{fact}\label{fact:threshComp}
Complement of a threshold graph is a threshold graph.
\end{fact}
\begin{definition}{Threshold dimension:}
A \emph{threshold cover} of a graph $G$ is a set of threshold graphs
$G_i$, $i=1,2,\ldots,k$ on the same vertex set as $G$ such that
$E(G)=E(G_1)\cup E(G_2)\cup \cdots \cup E(G_k)$.  The \emph{threshold
dimension} $t(G)$ is the least integer $k$ such that a threshold
cover of size $k$ exists. 
\end{definition}
Chv\'{a}tal and Hammer \cite{aggregationInequalityChvatal} introduced
threshold graphs and threshold dimension for their application in
set-packing problems.  In \cite{complexityPartialOrderDimnYannakakis},
Yannakakis showed that to determine if the threshold dimension of a
graph is at most 3 is NP-complete even for the class of split graphs.

For a graph $G$ let $G_i$, $1\le i\le k$ be graphs on the same vertex set
as $G$ such that $E(G)=E(G_1)\cap E(G_2)\cap\cdots\cap E(G_k)$. Then
we say that $G$ is the \emph{intersection} graph of $G_i$ s for $1\le
i\le k$ and denote it as $G=\bigcap_{i=1}^{k}G_i$. 
\begin{fact}
From Fact \ref{fact:threshComp} it is easy to see that threshold dimension of a graph $G$ is the smallest integer $k$ such that the complement graph $\overline{G}$ can be represented as the intersection of $k$ threshold graphs.
Also, if $G=\bigcap_{i=1}^{k}G_i$,
then $t(\comp{G})\le\sum_{i=1}^kt(\comp{G_i})$.
\end{fact}

\begin{lemma}\label{lem:threshSub}
Let $G$ be a split graph. Let $G'$ be a threshold supergraph of $G$. Then
we can construct another threshold graph $H$ such that
$G\subseteq H\subseteq G'$ and $\ind(H)=\ind(G)$.
\end{lemma}
\begin{proof}
First we observe that $\clq(G)\subseteq\clq(G')$. The graph $H$ is obtained
as follows: $\clq(H)=\clq(G)$ and $\ind(H)=\ind(G)$. For each
$u\in\ind(H)$, $N(u,H)=N(u,G')\cap\clq(G)$. By definition, $N(u,G)\subseteq
N(u,H)\subseteq N(u,G')$. Therefore $G\subseteq H\subseteq G'$. 

Now we will show that $H$ is a threshold graph. Suppose there exist
$u,v\in\ind(H)$, such that neither $N(u,H)\subseteq N(v,H)$ nor
$N(v,H)\subseteq N(u,H)$. There exist two vertices $u',v'\in\clq(H)$
such that $u'\in N(u,H)\setminus N(v,H)$ and $v'\in N(v,H)\setminus
N(u,H)$. This implies $u'\in N(u,G')\setminus N(v,G')$ and $v'\in
N(v,G')\setminus N(u,G')$, which in turn implies that $u'uvv'$ forms
an induced $P_4$ in $G'$. But, by Fact \ref{fact:tgDefn}, this is a
contradiction since $G'$ is a threshold graph.
\qed
\end{proof}

\subsection{Posets}
A \emph{partially ordered set} (or \emph{poset}) $P=(S,\leq_P)$ consists
of a non-empty finite set $S$ and a reflexive, antisymmetric and
transitive binary relation $\leq_P$ on $S$. $S$ is called the
\emph{ground set} of $P$. If $x\leq_P y$ or $y\leq_P x$ then $x$ and $y$
are said to be comparable. Otherwise we say that they are incomparable
and we denote this relation as $x||_Py$. We write $x<_P y$ when $x\leq_P
y$ and $x\neq y$. 

A \emph{totally ordered set} is a poset in which every two elements are
 comparable. A \emph{linear extension} $L$ of a poset $P$ is a totally
ordered set $(S,\leq_L)$ which satisfies: \(x\leq_Py \Longrightarrow
x\leq_Ly\). Let $L(u)=|\{v|v\le_Lu\}|$ denote the \emph{index} of the
element $u$ in the totally ordered set $L$.
 
A \emph{realizer} of a poset $P$ is a set of linear extensions of $P$,
say $\rlzr:L_1,L_2,\ldots,L_k$ which satisfy the following condition:
if $x||_Py$ then there exists two linear extensions $L_i,L_j\in\rlzr$
such that $x<_{L_i}y$ and $y<_{L_j}x$. The \emph{poset dimension} of $P$
denoted by $\dim(P)$ is the minimum integer $k$ such that there exists
a realizer of $P$ of cardinality $k$. It was introduced
by Dushnik and Miller \cite{partiallyOrderedSetsDushnik}. The poset
dimension problem is to decide for a given poset and integer $d$ whether
the dimension of the poset is at most $d$. For a
survey on dimension theory of posets see Trotter's
monograph \cite{combPosetsDimensionTrotter} or survey
paper \cite{partiallyOrderedSetsTrotter}. 

In \cite{complexityPartialOrderDimnYannakakis} Yannakakis studied the
complexity of the partial order dimension problem and its consequences on
various graph parameters. He proved that it is NP-complete to determine
whether the dimension of a partial order is at most 3. He then used
some simple reductions to extend this result to the problems of
determining the threshold dimension, boxicity and cubicity of graphs.
Recently in \cite{hardnessPODHegdeJain} Hegde and Jain reduced the fractional
chromatic number problem to the poset dimension problem to show that it is hard
to even approximate the dimension of a partial order. To state more
precisely,
\begin{theorem}{\textnormal{\cite{hardnessPODHegdeJain}}}\label{thm:jain}
 There exists no polynomial-time algorithm to approximate the poset dimension on an N-element set with a factor of $O(N^{0.5-\epsilon})$ for any $ \epsilon >0$, unless $NP=ZPP$.  
\end{theorem}

\subsection{Our Results}
In this paper we will show that 
\begin{enumerate}
 \item There exists no polynomial-time algorithm to approximate the threshold dimension of a graph on $n$ vertices with a factor of $O(n^{0.5-\epsilon})$ for any $ \epsilon >0$, unless $NP=ZPP$.
 \item There exists no polynomial-time algorithm to approximate the boxicity of a graph on $n$ vertices with a factor of $O(n^{0.5-\epsilon})$ for any $ \epsilon >0$, unless $NP=ZPP$.
 \item There exists no polynomial-time algorithm to approximate the cubicity of a graph on $n$ vertices with a factor of $O(n^{0.5-\epsilon})$ for any $ \epsilon >0$, unless $NP=ZPP$.
 \item If $G$ is a split graph then it is NP-complete to determine whether $\boxi(G) \leq 3$.
\end{enumerate}

\section{Preliminaries}
Let $G$ be a simple finite undirected graph on $n$ vertices. The vertex set of $G$ is denoted as $V(G)$ and the edge set of $G$ is denoted as $E(G)$. For each vertex $v \in V(G)$ let $N(v,G)$ denote the set
of vertices in $V(G)$ to which $v$ is adjacent. Whenever there is no
ambiguity regarding the graph under consideration, we will use the abbreviated
notation $N(v)$. A graph $H$ is said to be a subgraph of $G$ if and only
if $V(H)\subseteq V(G)$ and $E(H) \subseteq E(G)$. In this paper we will
use the notation $H \subseteq G$ to denote $H$ is a subgraph of $G$. Let
$V'\subseteq V$. $G[V']$ denotes the induced subgraph of $G$ on the vertex
set $V'$. For a positive integer $k$, let $[k]$ denote the set $\{1,2,\ldots,k\}$.

Suppose $I$ is an interval graph. Let us consider an interval representation of $I$. Without loss of generality we can assume that the endpoints of each interval are integers. For any vertex $u$, let $l(u)$ and $r(u)$ denote the integers corresponding to the left endpoint and right endpoint respectively of the interval corresponding to $u$.

\begin{property}{Helly property of intervals:}
Suppose $A_1,A_2,\ldots,A_k$ is a finite set of intervals on the real
line with pairwise non-empty intersection. Then there exists a common point
of intersection for all the intervals i.e. $\bigcap_{i=1}^k
A_i\ne\emptyset$.
\end{property}

Let $I_1,I_2,\cdots,I_k$ be $k$ interval graphs (unit interval graphs)
such that $G=\bigcap_{i=1}^k I_i$. Then $I_1,I_2,\cdots,I_k$ is called
an interval (unit interval) representation of $G$. Boxicity can be stated
in terms of intersection of interval graphs as follows:
\begin{lemma}{\textnormal{Roberts \cite{recentProgressesInCombRoberts}}}\label{lem:intbox}
The boxicity of a graph $G$ is the minimum positive integer $b$ such that $G$ can be represented as the intersection of $b$ interval graphs. Moreover, if $G=\bigcap_{i=1}^{m}G_i$ for some graphs $G_i$ then $\boxi(G)\leq\sum_{i=1}^{m}\boxi(G_i)$.  
\end{lemma}
Similarly cubicity can be stated in terms of intersection of unit interval graphs as follows:
\begin{lemma}{\textnormal{Roberts \cite{recentProgressesInCombRoberts}}}\label{lem:intcub}
The cubicity of a graph $G$ is the minimum positive integer $b$ such that $G$ is the intersection of $b$ unit interval graphs. Moreover, if $G=\bigcap_{i=1}^{m}G_i$ for some graphs $G_i$ then $\cub(G)\leq\sum_{i=1}^{m}\cub(G_i)$.  
\end{lemma}
The \emph{boxicity problem} is defined to be the problem of computing the boxicity for a given graph $G$.

%

%
%
%
%
%
%
\section{Characteristic Poset of a Split Graph}
In this section, we will introduce the concept of the characteristic
poset of a split graph and we will relate the threshold dimension and the
boxicity of split graphs to the dimension of this poset.
\begin{definition}\label{def:charPoset} 
Let $G$ be a split graph with $\ind(G)$ and $\clq(G)$
being the independent set and clique respectively. Let
$\gset(G)=\{N(u,G)|u\in\ind(G)\}$. The characteristic
poset of $G$ is $P=(\gset(G),\subseteq)$, i.e. the set of
neighborhoods of the independent set vertices ordered by inclusion.
\end{definition}
Note that the characteristic poset is unique to a split graph and 
by Fact \ref{fact:tgDefn}, we can infer that the characteristic poset
is a totally ordered set if and only if the split graph is a threshold graph.

\begin{theorem}\label{thm:charThresh}
Let $P$ be the characteristic poset of the split graph $G$. Then,
$\dim(P)\le t(\comp{G})$.
\end{theorem}
\begin{proof}
Let $t(\comp{G})=k$. Suppose $\thresh:T_1,T_2,\ldots,T_k$ is a set
of threshold graphs such that $\bigcap_{i=1}^kT_i=G$. From each
$T_i$, we will construct linear extension $L_i$ of $P$ such that
$L_1,L_2,\ldots,L_k$ form a realizer of $P$.

From Lemma \ref{lem:threshSub} we can assume that $\ind(T_i)=\ind(G)$
for $1 \le i \le k$. For each $T_i$ let $\gset(T_i)=\{N(u,T_i)|u\in
\ind(G)\}$. Consider the function $f_i:\gset(G)\longrightarrow\gset(T_i)$
where, for $X\in\gset(G)$, $f_i(X)$ is the smallest subset in
$\gset(T_i)$ containing $X$. Note that $f_i$ is well-defined: For each
$X\in\gset(G)$, there exists an $X'\in\gset(T_i)$ such that $X\subseteq
X'$ since $T_i$ is a supergraph of $G$. Moreover, the smallest subset
$f_i(X)$ is unique since $\gset(T_i)$ is a totally ordered set with
respect to set inclusion. We define
$L_i$ as follows: For any two distinct elements $X,Y\in\gset(G)$, 
\begin{enumerate}
\item If $f_i(X)\subset f_i(Y)$, then, $X<_{L_i}Y$.
\item If $f_i(X)=f_i(Y)$ and $X<_PY$, then, $X<_{L_i}Y$.
\item If $f_i(X)=f_i(Y)$ and $X||_PY$, then, we either make $X<_{L_i}Y$ or
$Y<_{L_i}X$.
\end{enumerate}
Since $T_i$ is a threshold graph, we observe that
\begin{eqnarray*}
X\subseteq Y &\Longrightarrow& f_i(X)\subseteq f_i(Y)\\
&\Longrightarrow& X\le_{L_i}Y
\end{eqnarray*}
Hence, $L_i$s are linear extensions of $P$. Suppose $X||_PY$, then there
exist $u,v\in\ind(G)$ such that $N(u,G)=X$ and $N(v,G)=Y$ and therefore there
exist $u',v'\in\clq(G)$ such that $u'\in N(u,G)\setminus
N(v,G)$ and $v'\in N(v,G)\setminus N(u,G)$. Since $\bigcap_{i=1}^kT_i=G$,
there exist two threshold graphs $T_j,T_l\in\thresh$ such that 
$u'\notin N(v,T_j)$ and $v'\notin N(u,T_l)$. This implies that
$f_j(Y)\subset f_j(X)$ and $f_l(X)\subset f_l(Y)$. Therefore, $Y<_{L_j}X$
and $X<_{L_l}Y$. Hence, we have proved that $L_i$s form a realizer of $P$.
\qed
\end{proof}
\begin{lemma}\label{lem:splitInt}
Let $G$ be a split graph. Let $G'$ be an interval supergraph of $G$. Then
we can construct a split interval graph $H$ such that
$G\subseteq H\subseteq G'$ and $\ind(H)=\ind(G)$.
\end{lemma}
\begin{proof}
Consider an interval representation of $G'$ such that it satisfies the
following two properties: (1) None of the intervals used is a single point
interval. (2) No two intervals share a common end point.
It is easy to see that such an interval representation can be constructed from
any given interval representation in polynomial time. Now let $x \in
\ind(G)$. Clearly $\{x\} \cup N(x,G)$ induces a clique in $G$ and therefore
in $G'$. Let $f'(v)$ denote the interval assigned to the vertex $v$
in the interval representation chosen for $G'$. By Helly property of the
intervals, $\bigcap_{v \in \{x\} \cup N(x,G)}f'(v) \neq \emptyset$. From
properties (1) and (2) we can easily infer that $\bigcap_{v \in \{x\} \cup
N(x,G)}f'(v)$ is not a single point interval. Now we define the interval
graph $H$ on the vertex set $V(G)$, by assigning the interval $f(v)$
to each vertex $v \in V(G)$, defined as follows 
\[
f(v)=\left\{
\begin{array}{ll}
f'(v) & \forall v \in \clq(G),\\
P(v) & \forall v \in \ind(G),
\end{array}\right.
\]
where $P(v)$ is a point in $\bigcap_{x \in \{v\} \cup N(v,G)}f'(x)$. Note
that since $\bigcap_{x \in \{v\} \cup N(v,G)}f'(x)$ is not a single
point we can assume that $P(v) \neq P(u)$ for all distinct $u,v
\in \ind(G)$. Also note that for each $v \in \ind(G)$, $N(v,G) \subseteq
N(v,H)$ by the construction. Since we have only changed the intervals
corresponding to the vertices in $\ind(G)$, we infer that $G \subseteq
H$. On the other hand $f'(v) \supseteq f(v)$ for all $v \in V(G)$ and
therefore $H \subseteq G'$, as required. Moreover it is easy to see that
$\ind(G)$ induces an independent set in $H$. Hence, $H$ is a split graph
with the same partition as $G$. Therefore, $H$ is a split interval graph. 
\qed 
\end{proof}

\begin{lemma}\label{lem:splitIntThresh}
If $G$ is a split interval graph, then $t(\comp{G})\le2$.
\end{lemma}
\begin{proof}
Let us consider an interval representation of $G$. We will construct
 two threshold graphs $G_1$ and $G_2$ as follows. Let $l=\min_{u\in
V(G)}l(u)$ and $r=\max_{u\in V(G)}r(u)$ be the leftmost and the rightmost
 points respectively, in the interval representation of $G$. Now, to
define $G_1$, we change the intervals corresponding to $u\in\clq(G)$ by
redefining their left end points: $l(u)=l$, $\forall u\in \clq(G)$. We do
not disturb the intervals corresponding to the vertices in $\ind(G)$. Now
we claim that $G_1$ is a threshold graph: Clearly $\ind(G)$ induces an
independent set in $G_1$ also. Therefore let $\ind(G_1)=\ind(G)$. Let
$u,v\in\ind(G_1)$. It is easy to see that $N(u,G_1)\supseteq N(v,G_1)$
if $l(u)\le l(v)$ and therefore, for every $u,v\in\ind(G_1)$, we have
either $N(u,G_1)\subseteq N(v,G_1)$ or $N(v,G_1)\subseteq N(u,G_1)$. 

Similarly, let $G_2$ be obtained by letting $r(u)=r$, $\forall u\in
\clq(G)$, while keeping other end points unchanged. Again by construction,
$G_2$ is a threshold graph. It is easy to see that $G_1\cap G_2=G$: By
construction, $G_1\supseteq G$ and $G_2\supseteq G$ and if $(u,v)\notin E(G)$,
it is clear that in $G_1$ or in $G_2$, the intervals corresponding to
$u$ and $v$ are disjoint. 
\qed
\end{proof}

\begin{lemma}\label{lem:splitBoxThresh}
If $G$ is a split graph, then $t(\comp{G})\le2\boxi(G)$.
\end{lemma}
\begin{proof}
Let $\boxi(G)=k$ and $G_1,G_2,\ldots,G_k$ be interval graphs on the same
vertex set as $G$ such that $\bigcap_{i=1}^kG_i=G$. By Lemma
\ref{lem:splitInt}, we can assume that all the $G_i$s are split interval
graphs. By Lemma \ref{lem:splitIntThresh}, corresponding to each $G_i$, we
can construct two threshold graphs $T_{2i-1}$ and $T_{2i}$ such that
$G_i=T_{2i-1}\cap T_{2i}$. Therefore, we have $2k$ threshold graphs whose
intersection gives $G$. Hence, proved.
\qed
\end{proof}
Combining the above Lemma and Theorem \ref{thm:charThresh}, we have:
\begin{theorem}\label{thm:charBox}
 Let $P=(S,\leq_P)$ be a characteristic poset of the split graph $G$. Then $\dim(P) \leq 2\boxi(G)$.
\end{theorem}

\begin{remark}
We observe that the constructions in Theorem \ref{thm:charThresh} and
Lemmas \ref{lem:splitInt}, \ref{lem:splitIntThresh} and
\ref{lem:splitBoxThresh} can be achieved in polynomial time.
\end{remark}

\section{Hardness of Approximation} \label{sec:lower_bound}
Given poset $P$, we will construct a split graph $G_P$ such that $P$
 is isomorphic to the characteristic poset of $G_P$. Consider a
poset $P=(S,\le_{P})$ where $|S|=n$. Let $g:[n]\longrightarrow S$
be a bijective map. For convenience, we will assume that $S$ and
$[n]$ are disjoint sets. We define a split graph $G_P$ as follows:
$V(G_P)=S\cup[n]$. $\clq(G_P)=[n]$ and $\ind(G_P)=S$. For any $u\in S$
and $v\in [n]$, $(u,v)\in E(G_P) \Longleftrightarrow g(v) \le_P u$. Thus
$g(N(u,G_P))=\{x \in S | x \leq_P u\}$. It is easy to see that $P$
is isomorphic to the characteristic poset of $G_P$.
\begin{theorem}\label{thm:threshLB}
$\dim(P)\ge t(\comp{G_P})$.
\end{theorem}
\begin{proof}
Let $\dim(P)=k$. Suppose $L_1,L_2,\ldots,L_k$ form a realizer of $P$. We
 will construct threshold graphs $G_i$ corresponding to each $L_i$ for
$1\le i\le k$ such that $\bigcap_{i=1}^kG_i=G_P$. The $G_i$s are defined as
follows: $V(G_i)=S\cup[n]$ with $\clq(G_i)=[n]$ and $\ind(G_i)=S$. For
any $u\in S$ and $v\in [n]$, $(u,v)\in E(G_i) \Longleftrightarrow
g(v) \le_{L_i} u$. $G_i$ is a threshold graph because $L_i$ (a
totally ordered set) is the characteristic poset of $G_i$. 

Now, we will show that if $(u,v)\in E(G_P)$ then $(u,v)\in E(G_i)$ $\forall
i\in[k]$. Since $\clq(G_i)=\clq(G_P)$, any $u,v\in\clq(G_i)$ are adjacent in
$G_i$. Suppose $u\in\ind(G_P)$ and $v\in\clq(G_P)$,
\begin{eqnarray*}
(u,v)\in E(G_P) &\Longrightarrow& g(v)\le_Pu\\
&\Longrightarrow& g(v)\le_{L_i}u, \forall i\in[k]\\
&\Longrightarrow& (u,v)\in E(G_i), \forall i\in[k]
\end{eqnarray*}
Hence, each $G_i$ is a supergraph of $G_P$. Next we will show that if
$(u,v)\notin E(G_P)$ then there exists $G_j$ such that $(u,v)\notin
E(G_j)$. If $(u,v)\notin E(G_P)$ then either $u<_Pg(v)$ or $u||_P
g(v)$. In either case, there exists an $L_j$ such that
$u<_{L_j}g(v)$. By definition of $G_j$, $(u,v)\notin E(G_j)$. Hence,
proved.
\qed
\end{proof}
Combining Theorems \ref{thm:charThresh} and \ref{thm:threshLB}, we have the following result.
\begin{corollary}\label{cor:dim}
 $\dim(P)=t(\overline{G_P})$.
\end{corollary}
Cozzens and Halsey \cite{relationshipTgSplitCozzensHalsey}
 proved that the boxicity of any graph $G(V,E)$ is not more
than the threshold dimension of it's complement $\comp{G}$,
i.e. $\boxi(G)\le t(\comp{G})$. Hence, 
\begin{corollary}\label{cor:boxPoset}
$\dim(P)\ge\boxi(G_P)$.
\end{corollary}
\begin{remark}
We note that the construction in Theorem \ref{thm:threshLB} can be achieved
in polynomial time.
\end{remark}

\begin{theorem}\label{thm:approx}
There exists no polynomial-time algorithm to approximate the threshold
dimension of a split graph on $n$ vertices with a factor of
$O(n^{0.5-\epsilon})$ for any $ \epsilon >0$ unless $NP=ZPP$.
\end{theorem}
\begin{proof}
Suppose there exists an algorithm to compute the boxicity of a split
graph on $n$ vertices with approximation factor $O(n^{0.5-\epsilon})$. As
 we have seen for any poset $P$ on $N$ elements we can construct a
split graph $G_P$ on $n=2N$ vertices such that $t(\comp{G_P})=\dim(P)$
 by Corollary $\ref{cor:dim}$. This immediately implies that $\dim(P)$
 can be approximated within factor $O(n^{0.5-\epsilon})$. But, from
Theorem \ref{thm:jain} we know that there exists no polynomial-time
algorithm to approximate the poset dimension problem with a factor
$O(n^{0.5-\epsilon})$ for any $\epsilon>0$ unless $NP=ZPP$, a
contradiction. 
\qed
\end{proof}

\begin{theorem}
There exists no polynomial-time algorithm to approximate the boxicity of a split graph on $n$ vertices with a factor of $O(n^{0.5-\epsilon})$ for any $ \epsilon >0$ unless $NP=ZPP$.
\end{theorem}
\begin{proof}
The proof is similar to that of Theorem \ref{thm:approx}. From
Theorem \ref{thm:charBox} and Corollary \ref{cor:boxPoset}, we have
$\boxi(G_P)\le \dim(P)\le 2\boxi(G_P)$. The rest follows from Theorem
\ref{thm:jain}.
\qed
\end{proof}

\begin{corollary}
There exists no polynomial-time algorithm to approximate the cubicity
of a split graph on $n$ vertices with a factor of $O(n^{0.5-\epsilon})$
for any $ \epsilon >0$ unless $NP=ZPP$.
\end{corollary}
\begin{proof}
In \cite{upperboundCubicityBoxicitySunilAshik} it is shown that for any graph $G$ on $n$ vertices,
$\cub(G)\le\boxi(G)\ceil{\log_2n}$. Since
any representation of $G$ as the intersection of cubes also serves as an
intersection of boxes, it follows that $\cub(G)\ge\boxi(G)$. Hence, given a
poset $P$ and the corresponding split graph $G_P$ as constructed in Section
\ref{sec:lower_bound}, we have $\cub(G_P)/\ceil{\log_2n}\le \dim(P)\le2\cub(G_P)$. The
rest follows as in Theorem \ref{thm:approx}.
\qed
\end{proof}
\section{NP-Completeness of Boxicity of Split Graph}
The following theorem was proved by Yannakakis in
\cite{complexityPartialOrderDimnYannakakis}.
\begin{theorem}\textnormal{\cite{complexityPartialOrderDimnYannakakis}}\label{lem:yan}
It is NP-complete to determine if a given split graph has threshold
dimension at most 3.
\end{theorem}
We will reduce the threshold dimension problem of split graphs to the
problem of computing boxicity of a split graph. 
Let $H$ be any split graph. Let $|V(H)|=n$. We will construct another
split graph $G'$ in polynomial time such that $\boxi(G')=t(H)$. A split
graph $G$ is said to be a complete split graph if for all $u \in \ind(G)$
and $v \in \clq(G)$, $(u,v) \in E(G)$. Note that a complete split graph
is also a threshold graph. If $H$ is a complete split graph then we
take $G'=H$ since $\boxi(H)=t(H)=1$. So for the rest of the proof we
will assume that $H$ is not a complete split graph. Let $G=\comp{H}$
and $G_1$, $G_2$ be copies of $G$. Let $V(G_1)=\clq(G_1)\cup\ind(G_1)$
and $V(G_2)=\clq(G_2)\cup\ind(G_2)$. $V(G')=V(G_1)\cup V(G_2)$ and
$E(G')=E(G_1)\cup E(G_2)\cup \{(u,v)|u\in \clq(G_1), v\in\clq(G_2)\}\cup
\{(u,v)|u\in\clq(G_1), v\in\ind(G_2)\}\cup \{(u,v)|u\in\clq(G_2),
 v\in\ind(G_1)\}$. Clearly, $G'$ is a split graph with
$\clq(G')=\clq(G_1)\cup \clq(G_2)$.
\subsection{$\boxi(G')\le t(H)$}
Let $t(H)=k$ and $T_1,T_2,\ldots,T_k$ be a set of threshold graphs such
that $\bigcap_{i=1}^kT_i=G$. Due to Lemma \ref{lem:threshSub}, we can
assume that $\ind(T_i)=\ind(G)$. Now we construct interval
graphs $H_i$ corresponding to each $T_i$ as follows: Let $T_i^1$ and
$T_i^2$ be two copies of $T_i$. We assume that $V(G_1)=V(T_i^1)$ and
$V(G_2)=V(T_i^2)$. Let $V(H_i)=V(G_1)\cup V(G_2)$. Let
$g_i:\ind(T_i^j)\longrightarrow [n]$, $j=1,2$, be a function
which assigns to each vertex in the independent set of $T_i^j$ a distinct
number satisfying: $u,v\in\ind(T_i^j)$, $N(u,T_i^j)\subset
N(v,T_i^j)\Longrightarrow g_i(u)> g_i(v)$. We define another function
$h_i:\clq(T_i^j)\longrightarrow [n]$, $j=1,2$, as: $\forall
u\in\clq(T_i^j)$
\[
h_i(u)=\left\{
   \begin{array}{ll}
   0, & \textrm{if } N(u,T_i^j)\cap \ind(T_i^j)=\varnothing,\\
   \ds \max_{v\in N(u,T_i^j)\cap \ind(T_i^j)}g_i(v), & \textrm{otherwise}.
   \end{array}\right.
\]
Each $u\in\ind(T^1_i)$ is associated with the single point interval
$[g_i(u),g_i(u)]$ and $u\in\clq(T^1_i)$ with interval $[-n,
h_i(u)]$. Each $u\in\ind(T^2_i)$ is associated with the single point
interval $[-g_i(u),-g_i(u)]$ and $u\in\clq(T^2_i)$ with interval $[
 -h_i(u),n]$. Now $H_i$ is defined to be the intersection graph of this family of intervals which corresponds to $V(G_1) \cup V(G_2)$.
\begin{remark}\label{remark:T}
 $\clq(T_i^j)=\clq(G_j)$ and $\ind(T_i^j)=\ind(G_j)$ for $1 \le i \le k$ and $j=1,2$.
\end{remark}
\begin{lemma}\label{lemma:split}
 $H_i$ is a split graph with $\clq(H_i)=\clq(G')$ and $\ind(H_i)=\ind(G')$ for $1 \le i \le k$. 
\end{lemma}
\begin{proof}
 In view of the construction of $H_i$ clearly, $0$ is a common point
for intervals corresponding to all vertices $u\in\clq(T^1_i)
\cup\clq(T^2_i)$. Also, by definition of $g_i$, it follows that 
intervals corresponding to all vertices $u\in\ind(T^1_i)\cup\ind(T^2_i)$
are mutually disjoint. Hence, $\clq(H_i)=\clq(G')$ and
$\ind(H_i)=\ind(G')$. Therefore, $H_i$ is a split graph.
\qed
\end{proof}
%
%

\begin{lemma}\label{lem:induced}
 $H_i[V(G_1)]=T_i^1$ and $H_i[V(G_2)]=T_i^2$ for $1 \le i \le k$.
\end{lemma}
\begin{proof}
 Clearly $H_i[V(G_1)]$ is a split graph with $\ind(H_i[V(G_1)])=\ind(T_i^1)$ and $\clq(H_i[V(G_1)])=\clq(T_i^1)$. By construction it is easy to see that $E(H_i[V(G_1)]) \supseteq E(T_i^1)$. Let $x \in \ind(T_i^1)$ and $y \in \clq(T_i^1)$ such that $(y,x) \notin E(T_i^1)$. Let $z \in \ind(T_i^1)$ be such that $(y,z) \in E(T_i^1)$. According to Fact \ref{fact:tgDefn} we have either $N(x,T_i^1) \subseteq N(z,T_i^1)$ or $N(x,T_i^1) \supseteq N(z,T_i^1)$. But since $y \notin N(x,T_i^1)$ and $y \in N(z,T_i^1)$ we can infer that $N(x,T_i^1) \subset N(z,T_i^1)$. It follows that $g_i(x) > g_i(z)$. Clearly $h_i(y) \le g_i(z) < g_i(x)$. Therefore $(x,y) \notin E(H_i[V(G_1)])$ and therefore $H_i[V(G_1)]=T_i^1$. A similar proof shows that $H_i[V(G_2)]=T_i^2$.
\qed
\end{proof}

\begin{lemma}\label{lower_bd}
 $\boxi(G')\le t(H)$. 
\end{lemma}

\begin{proof}
According to Lemma \ref{lemma:split}, $\clq(H_i)=\clq(G')$ and $\ind(H_i)=\ind(G')$ for $1 \le i \le k$.
Let $u\in\clq(G')$ and $v\in\ind(G')$. We consider the following cases:
\begin{enumerate}
\item $u\in\clq(G_1)$ and $v\in\ind(G_2)$: Then $(u,v)\in E(G')$ by
construction of $G'$. According to Remark \ref{remark:T} and by
construction of $H_i$, the interval corresponding to $u \in \clq(T_i^1)$
contains $[-n,0]$ and $v \in \ind(T_i^2)$ corresponds to a single point interval on the negative x-axis. It follows that $(u,v)\in
E(H_i)$ for $1 \le i \le k$.
\item $u\in\clq(G_2)$ and $v\in\ind(G_1)$: Similar to case 1.
\item $u\in\clq(G_1)$ and $v\in\ind(G_1)$: Note that $G'[V(G_1)]=G_1$ and by Lemma \ref{lem:induced},
$H_i[V(G_1)]=T_i^1$ for $1 \le i \le k$. Since $\bigcap_{i=1}^kT_i^1=G_1$ we have $\bigcap_{i=1}^kH_i[V(G_1)]=\bigcap_{i=1}^kT_i^1=G_1=G'[V(G_1)]$.
\item $u\in\clq(G_2)$ and $v\in\ind(G_2)$: Similar to case 3. We can show that $\bigcap_{i=1}^kH_i[V(G_2)]=\bigcap_{i=1}^kT_i^2=G_2=G'[V(G_2)]$.
\end{enumerate}
From the above points we can infer that if $(u,v) \in E(G')$ then $(u,v)
\in E(H_i)$ for $1 \le i \le k$ and if $(u,v) \notin E(G')$ then $(u,v)
\notin E(H_l)$ for some $l \in [k]$. Therefore $\bigcap_{i=1}^kH_i=G'$ and
hence $\boxi(G') \le k=t(H)$.
\qed
\end{proof}
\subsection{$\boxi(G')\ge t(H)$}
Let $\boxi(G')=l$ and $I_1, I_2,\ldots,I_l$ be interval graphs such that $\bigcap_{i=1}^lI_i=G'$. From
Lemma \ref{lem:splitInt} we can assume that each $I_i$ is a split
graph with $\ind(I_i)=\ind(G')$. Moreover,
\begin{remark}\label{remark:I_split}
 $I_i[V(G_1)]$ and $I_i[V(G_2)]$ are split graphs with $\ind(I_i[V(G_1)])=\ind(G_1)$ and $\ind(I_i[V(G_2)])=\ind(G_2)$ respectively for $1 \le i \le l$.
\end{remark}
We shall use the notation $T_C$ to denote a complete split graph.
\begin{lemma}\label{Lem:TD}
 With respect to an interval representation of $I_i$, let $u_l$ and $u_r$
be the vertices corresponding to the leftmost and rightmost intervals
respectively, among the vertices in $\ind(I_i)$.
\begin{enumerate}
 \item If $u_l\in \ind(G_1)$ and $u_r\in \ind(G_2)$ then $t(\comp{I_i[V(G_1)]})=1$ and $t(\comp{I_i[V(G_2)]})=1$. 
 \item If $u_l\in \ind(G_2)$ and $u_r\in \ind(G_1)$ then $t(\comp{I_i[V(G_1)]})=1$ and $t(\comp{I_i[V(G_2)]})=1$. 
 \item If $u_l,u_r\in \ind(G_1)$ then $t(\comp{I_i[V(G_1)]})\le2$ and $I_i[V(G_2)]=T_C$.
 \item If $u_l,u_r\in \ind(G_2)$ then $I_i[V(G_1)]=T_C$ and $t(\comp{I_i[V(G_2)]})\le2$.
\end{enumerate}
\end{lemma}
\noindent \textbf{Proof(1):} First we will prove that
$I_i[V(G_1)]$ is a threshold graph, which, by Fact \ref{fact:threshComp},
implies $t(\comp{I_i[V(G_1)]})=1$. By assumption $r(u)< r(u_r)$ for all
$u\in\ind(I_i)$, $u \neq u_r$. Since $\ind(G_1)\cup \ind(G_2)$ induces an
independent set in $I_i$ we have $r(u)<l(u_r)$ for all $u\in\ind(G_1)$
because otherwise $l(u_r) \leq r(u) < r(u_r)$ and hence intervals
corresponding to $u$ and $u_r$ intersect in the interval representation
of $I_i$. For any $v\in\clq(G_1)$, $r(v)\ge l(u_r)$ since by construction
of $G'$, $(v,u_r)\in E(G')$ and $G' \subseteq I_i$. Combining these
two observations, we get $r(u)<l(u_r)\le r(v)$ and thus $r(u)<r(v)$
for all $u\in\ind(G_1), v\in\clq(G_1)$. Suppose $u_1,u_2\in\ind(G_1)$
such that $r(u_1)\le r(u_2)$. Now for all $v\in\clq(G_1)$, $r(u_1)
\le r(u_2) <r(v)$. If $(u_1,v)\in E(I_i[V(G_1)])$ then $l(v)\le r(u_1)
\le r(u_2) <r(v)$. Hence $(u_2,v)\in E(I_i[V(G_1)])$ also. From this
and Remark \ref{remark:I_split}, it is clear that Fact \ref{fact:tgDefn}
holds for $I_i[V(G_1)]$. Therefore $I_i[V(G_1)]$ is a threshold graph.
Similarly, we can show that $t(\comp{I_i[V(G_2)]})=1$.

\noindent \textbf{Proof(2):} Similar to Proof of (1).

\noindent \textbf{Proof(3):} Since $\ind(G_1)\cup \ind(G_2)$ induces an
independent set in $I_i$, we have for all $u\in\ind(G_2)$, $l(u)>r(u_l)$
and $r(u)<l(u_r)$. Since by construction of $G'$ for all $v\in\clq(G_2)$,
 $(v,u_l)\in E(G')$, $(v,u_r)\in E(G')$ and $G' \subseteq I_i$, we have
 $l(v)\le r(u_l)$ and $r(v)\ge l(u_r)$. This implies $l(v)<l(u)\le
r(u)<r(v)$ for all $u\in\ind(G_2), v\in\clq(G_2)$. Hence all vertices
in $\ind(G_2)$ are adjacent to all vertices in $\clq(G_2)$. Now
$I_i[V(G_2)]$ is a complete split graph and hence $I_i[V(G_2)]=T_C$. On
the other hand by Remark \ref{remark:I_split}, $I_i[V(G_1)]$ is a
split interval graph. Hence from Lemma \ref{lem:splitIntThresh},
$t(\comp{I_i[V(G_1)]})\le 2$. \\
\noindent \textbf{Proof(4):} Similar to Proof of (3).
\begin{remark}\label{rem:compThresh}
Suppose $G$ is a split graph with $t(\comp{G})=k$. Let
 $\thresh:T_1,T_2,\ldots,T_k$ be a set of threshold graphs
such that $\bigcap_{i=1}^kT_i=G$. It is easy to see that there does not
exist a pair of graphs $T_i,T_j\in\thresh$ such that $T_i\subseteq
T_j$. Suppose this was not the case, then, $G=\bigcap_{l=1,l\ne j}^kT_l$,
i.e. we could discard $T_j$, thus contradicting the minimality of $k$.
\end{remark}
\begin{lemma}\label{upper_bd}
 $\boxi(G')\ge t(H)$.
\end{lemma}
\begin{proof}
Based on Lemma \ref{Lem:TD}, we can infer that $I_i[V(G_1)]$ belongs to
exactly one of the following 3 cases: 1) $t(\comp{I_i[V(G_1)]})=1$ and
$I_i[V(G_1)] \neq T_C$. 2) $t(\comp{I_i[V(G_1)]})\le 2$. 3)
$I_i[V(G_1)]=T_C$. Let $l_1,l_2,l_3$ be such that $l_j$ denotes the number
of times $I_i[V(G_1)]$ belongs to case $j$ for $1 \le i \le l$ and $1 \le j
\le 3$. Clearly $l_1+l_2+l_3=l$. Recall that $H$ is not a complete split graph. Therefore there exists some $i \in [l]$ such that $I_i \neq T_C$. Note that $G_1=\bigcap_{i=1}^l I_i[V(G_1)] $ and therefore $t(\comp{G_1}) \le \sum_{i=1}^l t(\comp{I_i[V(G_1)]})$ $\le l_1+2l_2+l_3 t(\comp{T_C})$. Since any threshold graph $T$ which is a supergraph of $\comp{H}$ is a subgraph of $T_C$, by Remark \ref{rem:compThresh}, $T_C$ can be discarded and therefore, we can ignore the term $l_3 t(\comp{T_C})$ in the above expression. Hence we get $t(\comp{G_1}) \le l_1+2l_2$.\\
We can get 3 similar cases for $I_i[V(G_2)]$. Let $l_j'$ denotes the number of times $I_i[V(G_2)]$ belongs to case $j$ for $1 \le i \le l$ and $1 \le j \le 3$. Clearly $l_1'+l_2'+l_3'=l$. From Lemma \ref{Lem:TD}, it is easy to see that $l_1'=l_1$, $l_2'=l_3$ and $l_3'=l_2$. Therefore $t(\comp{G_2}) \le \sum_{i=1}^l t(\comp{I_i[V(G_2)]})$ $\le l_1+2l_3$. Hence realizing that $\comp{G_1}$ and $\comp{G_2}$ are isomorphic to $H$,
\begin{eqnarray*}
2t(H)=t(\comp{G_1})+t(\comp{G_2})\le 2(l_1+l_2+l_3)=2l.
\end{eqnarray*}
Hence, we get $t(H)\le l=\boxi(G')$. 
\qed
\end{proof}

\begin{theorem}
It is NP-complete to determine if a given split graph has boxicity at most 3.
\end{theorem}
\begin{proof}
We reduce the problem of determining the threshold dimension of a
split graph to this problem. Given a split graph $H$ we can construct another split graph $G'$ in polynomial time such that $\boxi(G')=t(H)$ by  Lemma \ref{lower_bd} and Lemma \ref{upper_bd}. The rest follows from Theorem \ref{lem:yan}.
\qed
\end{proof}


\begin{thebibliography}{10}
\expandafter\ifx\csname url\endcsname\relax
  \def\url#1{\texttt{#1}}\fi
\expandafter\ifx\csname urlprefix\endcsname\relax\def\urlprefix{URL }\fi

\bibitem{upperboundCubicityBoxicitySunilAshik}
L.~S. Chandran, K.~A. Mathew, An upper bound for cubicity in terms of boxicity,
  Disc. Math. 309~(8) (2009) 2571--2574.

\bibitem{aggregationInequalityChvatal}
V.~Chv\'{a}tal, P.~L. Hammer, Aggregation of inequalities in integer
  programming, in: Ann. Discrete Math, 1977.

\bibitem{phdThesisCozzens}
M.~B. Cozzens, Higher and multi-dimensional analogues of interval graphs, Ph.D.
  thesis, Department of Mathematics, Rutgers University, New Brunswick, NJ
  (1981).

\bibitem{relationshipTgSplitCozzensHalsey}
M.~B. Cozzens, M.~D. Halsey, The relationship between the threshold dimension
  of split graphs and various dimensional parameters, Disc. Appl. Math. 30
  (1991) 125--135.

\bibitem{partiallyOrderedSetsDushnik}
B.~Dushnik, E.~W. Miller, Partially ordered sets, Amer. J. Math 6~(3) (1941)
  600--610.

\bibitem{splitGraphFoldes}
S.~Foldes, P.~L. Hammer, Split graphs, in: Proceedings of the 8th South-Eastern
  Conference on Combinatorics, Graph Theory and Computing, 1977.

\bibitem{algGraphTheoryPerfectGraphsGolumbic}
M.~C. Golumbic, Algorithmic Graph Theory and Perfect Graphs, Academic Press,
  New York, 1980.

\bibitem{hardnessPODHegdeJain}
R.~Hegde, K.~Jain, The hardness of approximating poset dimension, Electronic
  Notes on Discrete Mathematics 29 (2007) 435--443.

\bibitem{specialPlanarSatisfiabilityProbNPKratochvil}
J.~Kratochv{\'{\i}}l, A special planar satisfiability problem and a consequence
  of its {N}{P}-completeness, Disc. Appl. Math. 52 (1994) 233--252.

\bibitem{recentProgressesInCombRoberts}
F.~S. Roberts, Recent Progresses in Combinatorics, chap. On the boxicity and
  cubicity of a graph, Academic Press, New York, 1969, pp. 301--310.

\bibitem{combPosetsDimensionTrotter}
W.~T. Trotter, Combinatorics and partially ordered sets: Dimension Theory, The
  Johns Hopkins University Press, Baltimore, Maryland, 1992.

\bibitem{partiallyOrderedSetsTrotter}
W.~T. Trotter, Graphs and partially ordered sets: recent results and new
  directions, in: Surveys in graph theory (San Fransisco, CA, 1995), Congr.
  Numer. 116, 1996.

\bibitem{canonicalPartitionTyshkevich}
R.~I. Tyshkevich, A.~A. Chernyak, Canonical partition of a graph defined by the
  degrees of its vertices, in: (In Russian) Isv. Akad. Nauk BSSR, Ser.
  Fiz.-Mat. Nauk, 1979.

\bibitem{complexityPartialOrderDimnYannakakis}
M.~Yannakakis, The complexity of the partial order dimension problem,
  SIAM~J.~Alg.~Disc.~Math. 3~(3) (1982) 351--358.

\end{thebibliography}
%
\end{document}